\documentclass[12pt,leqno]{article}
\usepackage{graphicx,psfrag,amssymb,amsmath,amsthm,makeidx}
\usepackage{epstopdf}

\def\IR{{\mathbb R}}
\def\Tor{\text{\rm Tor}}
\def\pQ{{\partial Q}}
\def\cH{\mathcal{H}}
\def\cM{\mathcal{M}}
\def\cS{\mathcal{S}}
\def\cT{\mathcal{T}}
\def\bm{\mathbf{m}}
\def\eps{\varepsilon}
\def\la{\langle}
\def\ra{\rangle}
\def\beq{\begin{equation}}
\def\eeq{\end{equation}}

\newtheorem{theorem}{Theorem}
\newtheorem{lemma}{Lemma}

\begin{document}

\title{Upgrading the Local Ergodic Theorem \\
for planar semi-dispersing billiards}
\author{ N. Chernov$^1$ and N. Sim\'anyi$^1$}

\maketitle

\begin{abstract}
The Local Ergodic Theorem (also known as the `Fundamental Theorem') gives
sufficient conditions under which a phase point has an open
neighborhood that belongs (mod 0) to one ergodic component. This
theorem is a key ingredient of many proofs of ergodicity for
billiards and, more generally, for smooth hyperbolic maps with
singularities. However, the proof of that theorem relies upon a
delicate assumption (Chernov-Sinai Ansatz), which is difficult to
check for some physically relevant models, including gases of hard
balls. Here we give a proof of the Local Ergodic Theorem for two
dimensional billiards without using the Ansatz.
\end{abstract}


\footnotetext[1]{Department of Mathematics, University of Alabama at
Birmingham, AL, 35294, USA; Email:$\ $ chernov@math.uab.edu and
simanyi@math.uab.edu.}

\section{Introduction}

In this work we make a step toward a complete solution (yet to be
achieved) of the celebrated Boltzmann-Sinai ergodic hypothesis. The
latter asserts \cite{Si63} that every system of $n\geq 2$ hard balls
on a torus of dimension $d\geq 2$ is ergodic (provided the trivial
first integrals are eliminated). This model reduces to the motion of
a billiard particle in a $d(n-1)$-dimensional torus bouncing off
$n(n-1)/2$ cylindrical obstacles (the billiard particle hits a
cylinder whenever two balls collide). Billiards with cylindrical
walls belong to a more general category of semi-dispersing
billiards, where a particle moves in a container with concave (but
not necessarily strictly concave) boundaries.

We remark that in the case $n=2$ the cylinders actually become
spheres, i.e.\ any system of $2$ hard balls reduces to a billiard
particle in a torus with a spherical obstacle. Such billiards belong
to a more special class of dispersing billiards, where a particle
moves in a container with strictly concave walls.

Dispersing billiards are always completely hyperbolic and ergodic
\cite{SC87}, but for semi-dispersing billiards this may not be true.
For example, a billiard in a 3-torus with a single cylindrical wall
has zero Lyapunov exponents and is not ergodic; on the other hand, 2
transversal cylindrical walls within a 3-torus ensure hyperbolicity
and ergodicity \cite{KSS89}. For the systems of $n \geq 3$ hard
balls, one has to carefully explore the geometry of the cylindrical
walls in order to derive hyperbolicity and ergodicity.

There are two complications in the study of hard balls (or more
generally, semi-dispersing billiards). One is caused by the
\emph{singularities} of the dynamics -- these happen during
simultaneous multiple collisions of $\geq 3$ balls and during
grazing (tangential) collisions. In the phase space, singular points
make submanifolds of codimension one. The other complication is
caused by \emph{non-hyperbolicity} (i.e.\ the existence of zero
Lyapunov exponents) at some phase points. Such points make various
structures, ranging from  smooth submanifolds to Cantor-like subsets
of the phase space.

Powerful techniques have been developed to handle these two
complications separately (singularities and non-hyperbolicity), but
the combination of the two still presents an unmanageable situation.
More precisely, if non-hyperbolic sets and singularities intersect
in a subset of positive $[2d(n-1)-2]$-dimensional measure, then
modern proofs of ergodicity stall. On the other hand, such
substantial overlaps between singularities and non-hyperbolic sets
appear very unlikely (physically); they are regarded as
`conspiracy'.

To bypass this scenario in an early work, Ya.~Sinai and N.~Chernov
\cite{SC87} \emph{assumed} that almost every point on the
singularity manifolds (with respect to the intrinsic Lebesgue
measure) was completely hyperbolic. Under this assumption (now
referred to as Chernov-Sinai Ansatz) they proved the so-called Local
Ergodic Theorem (also called `Fundamental Theorem'), which later
became instrumental in the proofs of ergodicity for various
billiards \cite{BCST,KSS90,LW}. It gives sufficient (and easily
verifiable) conditions under which a phase point has an open
neighborhood which belongs (mod 0) to one ergodic component.

A.~Kr\'amli, N.~Sim\'anyi and D.~Sz\'asz built upon the results of
\cite{SC87} and established the ergodicity for systems of $n=3$ hard
balls in any dimension \cite{KSS91} and for $n=4$ hard balls in
dimension $d\geq 3$ \cite{KSS92}; in particular they verified
Chernov-Sinai Ansatz in these cases. However, their techniques could
not be extended to $n \geq 5$. The situation called for novel
approaches.

A partial breakthrough was made by Sim\'anyi and Sz\'asz when they
invoked ideas of algebraic geometry to rule out various
`conspiracies' (at least for generic systems of hard balls), which
were in the way of proving hyperbolicity and ergodicity. Precisely,
they assumed that the balls had arbitrary masses $m_1,\dots,m_n$
(but the same radius $r$) and proved \cite{SS99} complete
hyperbolicity at a.e.\ phase point for generic vectors of `external
parameters' $(m_1, \ldots, m_n, r)$; the latter needed to avoid some
exceptional submanifolds of codimension one in $\IR^{n+1}$, which
remained unspecified and unknown. Later Sim\'any used
\cite{Sim03,Sim04} the same approach to prove Chernov-Sinai Ansatz
and ergodicity for generic systems of hard balls (in the above
sense). He also established hyperbolicity for systems of hard balls
of arbitrary masses \cite{Sim02}.

Thus the Boltzmann-Sinai ergodic hypothesis is now proved for
typical, or generic, systems of hard balls. This seems to be a
comforting settlement in both topological and measure-theoretic
senses, but it falls short of solving physically relevant problems,
as there is no way to check whether any particular system of hard
balls is ergodic or not. Most notably, for the system of balls with
all equal masses (which lies in the foundation of statistical
mechanics) the ergodicity remains open.

In an attempt to extend his results to ALL gases of hard balls
(without exceptions), Simanyi developed \cite{Sim07} a new approach
based on purely dynamical (rather than algebro-geometric) ideas;
this allowed him to derive ergodicity from Chernov-Sinai Ansatz for
all hard ball systems. Thus the Boltzmann-Sinai hypothesis is now
solved conditionally, modulo the Ansatz. It remains to prove Ansatz,
or alternatively, derive Local Ergodic Theorem without Ansatz. (As a
side remark, it is ironic that Ansatz, which originally seemed to be
just a convenient and temporary technical assumption, now remains
the only unresolved issue in the whole picture.)

Here we make another step toward a final solution of the classical
ergodic hypothesis: we derive Local Ergodic Theorem without Ansatz
for arbitrary semi-dispersing billiards in dimension two. Our method
does not yet apply to higher dimensions, but we are working on this.

\section{Statement of the result}

A planar (two-dimensional) billiard is a dynamical system where a
point $q$ moves freely with unit velocity $v$, $\|v\|=1$, in a
bounded connected domain $Q \subset \IR^2$ or $Q \subset \Tor^2$ and
reflects off its boundary $\pQ$ by the rule
\beq \label{v+v-}
  v^+ = v^- - 2\, \la n(q),v^-\ra\, n(q)
\eeq
where $v^+$ and $v^-$ denote the postcollisional and precollisional
velocities, $n(q)$ is the inward unit normal vector to $\pQ$ at the
collision point $q \in \pQ$, and $\la \cdot, \cdot \ra$ is the
scalar product in $\IR^2$.

As usual, $\pQ$ is a finite union of $C^3$ compact curves that can
only intersect at common endpoints (which make corners of the table
$Q$). Whenever the particle hits a corner point $q \in \pQ$, there
are two normal vectors to $\pQ$, thus the rule \eqref{v+v-} gives
two possible continuations (two branches) of the billiard
trajectory. Of course, this is an exceptional event (a singularity,
see below).

A billiard table $Q$ is semi-dispersing if every smooth component of
$\pQ$ is convex (but not necessarily strictly convex) inward. We
also suppose that the set of inflection points $q \in \pQ$ (where
the curvature of $\pQ$ vanishes) is a finite union of straight line
segments (flat sides of $Q$) and some isolated points. A simple
example is a polygon with one or several convex ovals removed from
its interior. In semi-dispersing billiards, collisions cannot
accumulate \cite{G,Va}, i.e.\ within any finite time period the
particle experiences finitely many collisions, hence its trajectory
is always well defined (though it might be multiply defined, due to
corner points).

The phase space of the billiard system is a compact three
dimensional manifold $\Omega = Q \times S^1$, and the billiard flow
$\Phi^t \colon \Omega \to \Omega$ preserves a uniform measure $\mu$
on $\Omega$. The collision space
$$
   \cM=\{(q,v)
   \in \Omega \colon q\in \pQ,\,
   \la v,n(q) \ra \geq 0\}
$$
consists of all postcollisional velocity vectors at reflection
points. We define the first collision time $\tau(x) = \min
\{t>0\colon \Phi^t (x) \in \cM\}$ and the (first) collision map $T
(x) = \Phi^{\tau(x)+0} (x)$ that maps $\Omega$ onto $\cM$; its
restriction to $\cM$ is called the billiard map (or collision map).
Canonical coordinates on $\cM$ are $r$ and $\varphi$, where $r$ is
the arc length parameter on $\pQ$ and $\varphi\in [-\pi/2,\pi/2]$ is
the angle between $v$ and $n(q)$. The map $T \colon \cM \to \cM$
preserves the smooth measure $d\nu = \cos \varphi \, dr \,
d\varphi$.

For every $x=(q,v) \in \Omega$ we put $-x=(q,-v)$; similarly for
every $x=(q,v^+) \in \cM$ we put $-x =(q,-v^-)$, where $v^+$ and
$v^-$ are related by (\ref{v+v-}).

If $Q \subset \Tor^2$, we may have an unpleasant case of `infinite
horizon', where $\sup_x \tau(x) = \infty$. In that case we enlarge
$\cM$ to make the horizon finite \cite{SC87}.
Suppose $\Tor^2$ is obtained by
identifying the opposite sides of the boundary of a rectangle $K$;
then we add the set $\partial K \times S^1$ to $\cM$. In other
words, every time the particle crosses $\partial K$, we record a
`collision', though the particle keeps moving straight with the same
velocity (we call $\partial K$ a transparent wall). Now it is clear
that $\sup_x \tau(x) < \infty$.

The billiard flow $\Phi^t$ is a suspension flow over the base map $T
\colon \cM \to \cM$ under the ceiling function $\tau$; it is ergodic
if and only if $T$ is.

The flow $\Phi^t$ and the map $T$ are singular (non-differentiable)
whenever the particle hits a corner of $Q$ or makes a grazing
(tangential) collision with $\pQ$, i.e.\ whenever the next collision
point belongs to
${\cal S}_0 = \{(q,v)\in \Omega:\; \langle v,n(q)\rangle
= 0 \text{ or } q\in\Gamma^\ast\}$,
where $\Gamma^\ast$ denotes the set of corner
points (observe that $\varphi = \pm \pi/2$ at grazing collisions).
The singularity set $\cS_1 = T^{-1} (\cS_0)$ of $T$ is a finite
union of smooth compact curves in $\cM$ (it is exactly to ensure its
finiteness why we added the transparent wall to $\cM$). Similarly,
for each $n \neq 0$ the singularity set $\cS_n = T^{-n} (\cS_0)$
(which is part of the singularity set of the iterate $T^n$)
is a finite union of smooth compact curves in $\cM$.

In semi-dispersing billiards, the set $\cS_n$ consists of increasing
curves for $n<0$ and of decreasing curves for $n>0$; thus
singularity curves always intersect each other transversally at some
time in their lives. Points
$x \in \cM$ whose trajectories are singular in both future and past
(the so called `double-singularities') make a countable set, which
can be easily neglected in the studies of ergodic properties of $T$.
Accordingly, the singularities of the flow $\Phi^t$ are a countable
union of hypersurfaces in $\Omega$, and future singularities
intersect past singularities transversally.

Next we describe hyperbolic properties of $\Phi^t$ and $T$. A local
orthogonal manifold (LOM), also called wave front, denoted by
$\Sigma \subset \Omega$, is a smooth oriented curve $\gamma \subset
Q$ equipped with a family of unit normal vectors (note that there
are exactly two such families). The $\Phi^t$-image of a LOM is a
finite union of LOMs (sometimes having common endpoints) in
$\Omega$.

If the map $T$ is smooth on $\Sigma \subset \Omega$ then, slightly
abusing notation, we call $\Sigma^c = T(\Sigma) \in \cM$ a LOM as
well. Given a LOM $\Sigma^c \subset \cM$, we call $\Sigma \subset
\Omega$ the corresponding flow-sync LOM (the latter is not unique of
course).

We distinguish divergent, convergent, and flat LOMs, as determined
by the curvature of its carrier $\gamma \subset Q$. In
semi-dispersing billiards, future images of divergent LOMs are
always divergent and their sizes keep growing in time; this is the
cause of hyperbolicity. On the other hand, images of flat LOMs
remain flat as long as they collide with flat sides of $Q$; but they
become divergent immediately after a collision with a curved side of
$Q$.

We assume that $\pQ$ has non-zero curvature in at least one point;
otherwise $Q$ is a polygon and there are no hyperbolic points.
Billiards in generic polygons are ergodic \cite{KMS} (though it is
hard to construct explicit examples \cite{Vo}), but they are never
hyperbolic.

For $x \in \Omega$ and $a<b$, a trajectory segment $\Phi^{[a,b]}
(x)$ of the point $x$ is said to be sufficient if there is a
collision at some time $a<t<b$ with a curved side of $Q$ (the
curvature of $\pQ$ must be different from zero at the collision
point); if the segment $\Phi^{[a,b]} (x)$ passes through singular
points and branches out, then every branch must hit a curved side of
$Q$. A point $x \in \Omega$ is sufficient in the future (past) if
its semitrajectory $\Phi^{[0,\infty)} (x)$ (resp.,
$\Phi^{(-\infty,0]} (x)$) is sufficient. If a nonsingular point $x
\in \cM$ is sufficient (future or past), then in a vicinity $U_x$ of
$x$ almost every point $y \in U_x$ is hyperbolic (this follows from
the Poincar\'e theorem); thus sufficiency guarantees (local)
hyperbolicity.

\medskip \noindent \textbf{Chernov-Sinai Ansatz}. Almost every point
$x \in \cS_1$ (with respect to the one-dimensional Lebesgue measure
on $\cS_1$) is past sufficient or, equivalently, almost every point
$x \in \cS_{-1}$ is future sufficient. \medskip

Here is our main result:

\begin{theorem} \label{TmMain}
Let $x_0 \in \cM$ be a point whose entire trajectory
$\Phi^{(-\infty, \infty)} (x)$ passes through at most one
singularity and is sufficient. Then there exists an open
neighborhood $U_0$ of $x_0$ that belongs (mod 0) to one ergodic
component of the map $T$.
\end{theorem}

\bigskip

We note that the base neighborhood $U_0=U_0(x_0)$ in this theorem is
any open neighborhood $U_0$ of $x_0$ for which

\medskip

(i) $U_0$ is a subset of the neighborhood $U_{\eps_1}(x_0)$ of $x_0$
featuring Theorem 3.6 in \cite{KSS90}, where $0<\eps_1<1$ is any fixed
number, and

\medskip

(ii) $U_0$ admits a family

$$
\mathcal{G}^\delta=\left\{G_i^\delta \big|\; i=1,2,\dots,I(\delta)
\right\} \quad (0<\delta<\delta_0)
$$
of regular coverings with a small enough threshold $\delta_0>0$ as
explained in \cite{KSS90}.

\bigskip

All the existing proofs of the Local Ergodic Theorem \cite{SC87,KSS90}
assume the Ansatz, and we relax that assumption.

Given a particular 2d semi-dispersing billiard, one can verify its
ergodicity by showing that the set of sufficient points is connected
and has full measure. We note, however, that it is unknown if every
semi-dispersing billiard (excluding polygons) is ergodic (or even completely
hyperbolic), and our result will not solve this open problem,
because we cannot yet control the measure of insufficient points.
Proving that a.e.\ phase point in any semi-dispersing billiard is
sufficient amounts to showing that in every polygonal billiard a.e.\
trajectory is dense, but this is an old (and notoriously hard) open
problem.

An interesting result in this direction was obtained in \cite{CT}:
it was shown that billiards in any polygon where a `bump' or a
pocket is attached at every vertex are hyperbolic and ergodic. But
in that case the verification of Ansatz was trivial, as every
non-sufficient trajectory was periodic.

\section{Proof of the result}

We begin with a helpful geometric fact that gives a sufficient
condition under which two nearby phase points (points in $Q$
equipped with unit velocity vectors) belong to one divergent local
orthogonal manifold.

\begin{lemma} \label{Lm2.1}
Let $(q_1,v_1),\,(q_2,v_2)\in \IR^2 \times \IR^2$, $\|v_i\|=1$,
$\|q_1-q_2\|<\eps_0$, $\|v_1-v_2\|<\eps_0$, $\langle
q_1-q_2,\,v_1-v_2\rangle\ge 0$, with some fixed constant $\eps_0\ll
1$.  We claim that there are reals $\tau_1,\,\tau_2 \in \IR$,
$|\tau_i|<10000\eps_0$, such that the phase points
$(q_1+\tau_1v_1,\,v_1)$ and $(q_2+\tau_2v_2,\,v_2)$ can be included in
a divergent LOM $\Sigma \subset \IR^2 \times S^1$.
\end{lemma}

\emph{Proof}. We assume the strict inequality
$\langle q_1-q_2,\,v_1-v_2\rangle>0$. The general result then follows
by simply passing to the limit.

Let $O$ be the point of intersection of the lines
$l_1=\left\{q_1+tv_1\big|\; t\in\mathbb{R}\right\}$ and
$l_2=\left\{q_2+tv_2\big|\; t\in\mathbb{R}\right\}$. We may and shall
assume that $O\in\mathbb{R}^2$ is the origin of the plane $\mathbb{R}^2$.
Let $q_1=t_1v_1$, $q_2=t_2v_2$. The assumed inequality says that
$\langle q_1-q_2,\,v_1-v_2\rangle=(t_1+t_2)(1-\langle
v_1,v_2\rangle)>0$, thus $t_1+t_2>0$, so by symmetry we may assume
that $t_1>0$. We distinguish between two cases:

\textbf{Case 1}: $t_1\ge 5000\eps_0$. We take
$q_2+\tau_2v_2=t_1v_2$, i. e. $\tau_1=0$ and $\tau_2=t_1-t_2$.
Clearly, both $(q_1,v_1)$ and $(q_2+\tau_2v_2,v_2)$ are elements of
the (outer unit normal field of) the circle $\Sigma$ defined by the
equation $\|x\|=t_1$.

\textbf{Case 2}: $0<t_1<5000\eps_0$. We take
$q_1+\tau_1v_1=5000\eps_0v_1$, $q_2+\tau_2v_2=5000\eps_0v_2$, i.
e. $\tau_1=5000\eps_0-t_1$, $\tau_2=5000\eps_0-t_2$. The (unit
normal bundle of the) circle $\Sigma$ containing
$(q_1+\tau_1v_1,\,v_1)$ and $(q_2+\tau_2v_2,\,v_2)$ is now defined
by the equation $\|x\|=5000\eps_0$ in $\mathbb{R}^2$. \qed \medskip

Next we turn to the proof of Local Ergodic Theorem (without using
Ansatz). We follow the lines and notation of \cite{KSS90} that
presents one of the clearest and most complete proofs of that
theorem. For the given sufficient point $x_0$ we consider a small
enough open neighborhood $U_0=U_0(x_0)$ of $x_0$, as described right
after Theorem 1.

Given a divergent LOM $\Sigma \subset \Omega$ with a carrier $\gamma
\subset Q$, we use the metric on it generated by the distance along
the curve $\gamma$, and denote by $\|\cdot\|$ the corresponding norm
in its tangent space $\cT \Sigma$. For LOMs $\Sigma \subset \cM$ we
use the norm and metric on the corresponding flow-sync LOM's in
$\Omega$ (constructed right at the given point $x \in \Sigma$).

For any $x \in \Sigma \subset \cM$ denote by $D^n_{x,\Sigma}$ the
Jacobian of the map $T^n$ restricted to $\Sigma$ at $x$, in the
above norm. If $\Sigma$ is a divergent LOM (in our terminology the
word ``divergent'' always means ``not necessarily strictly divergent'',
and a similar convention applies to convergent LOMs), then
$D^n_{x,\Sigma} \ge 1$ for every $n\geq 1$. Denote
$$
  \kappa_{n,0} (x) = \inf_{\Sigma} D^n_{-T^n x,\Sigma}
$$
where the infimum is taken over all divergent LOMs through $-T^nx$;
this quantity is the minimal expansion of divergent LOMs on their
way from $-T^nx$ back to $-x$ (we note that the inf is actually
attained at the flat LOM, cf.\ \cite{CM}). Given $\delta>0$ we
denote
$$
  \kappa_{n,\delta} (x) = \inf_{\Sigma} \inf_{y\in\Sigma}
  D^n_{y,\Sigma}
$$
where the infimum is taken over all divergent LOMs $\Sigma$
through $-T^nx$
such that $T^n$ is smooth on $\Sigma$ and dist$(-x,
\partial T^n\Sigma) \leq \delta$ (of course, for any LOM $\Sigma$, the
boundary $\partial \Sigma$ consists of the two endpoints of that
LOM). We observe that $1 \leq \kappa_{n,\delta} (x) \leq
\kappa_{n,0} (x)$, and both $\kappa_{n,0} (x)$ and
$\kappa_{n,\delta} (x)$ are non-decreasing functions of $n$.

For $x\in \cM$ we denote
$$
   z_{\rm tub}(x) = \sup_{\Sigma} \{{\rm dist}(x,\partial \Sigma)\colon
   T\text{ is smooth on }\Sigma\}
$$
where the supremum is taken over all flat LOMs
$\Sigma\subset\text{int}\mathcal M$ through $x$; this is
the so-called radius of the maximal tubular neighborhood of the
billiard link joining $x$ with $Tx$.

Note that $z_{\text{tub}}(-Tx)=z_{\text{tub}}(x)$.

We denote by $\Sigma^u (x)$ and $\Sigma^s (x)$ the unstable and
stable manifolds through $x$; the former is a divergent LOM and the
latter a convergent one. We also put $r^\alpha (x) =\,$dist$(x,
\partial \Sigma^\alpha (x))$ for $\alpha = u,s$. It is known
\cite[Lemma~5.4]{KSS90} that for every semi-dispersing billiard
table $Q$ there exists a constant $c_3>0$ (using the notation of
\cite{KSS90}) such that if
$$
   x \in U^{\rm g} = U^{\rm g} (\delta) =
   \{y\in\cM\colon \ \ \forall n>0 \ \ z_{\rm tub}(-T^ny)
   \geq (\kappa_{n,c_3\delta}(y))^{-1}c_3\delta\}
$$
then $r^s(x) \geq c_3\delta$. Thus the points of $U^{\rm g}$ (`good set') have
stable manifolds of order $\delta$. A similar property holds for
unstable manifolds. The set of points with shorter stable manifolds
(`bad set') must be carefully analyzed. We put
\begin{align*}
   U^{\rm b} &= U^{\rm b} (\delta) = U_0 \setminus U^{\rm g} = \cup_{n
     \geq 1} U^{\rm b}_n\\
   U_n^{\rm b} &= U^{\rm b}_n (\delta) = \{y\in U^{\rm b} \colon z_{\rm tub}(-T^ny) <
   (\kappa_{n,c_3\delta}(y))^{-1} c_3\delta\}.
\end{align*}
A crucial fact in the proof of Local Ergodic Theorem is the
following tail bound: for any function $F(\delta) \to \infty$ as
$\delta \to 0$ the set
$$
   U^{\rm b}_{\omega} = U^{\rm b}_{\omega} (\delta) =
   \cup_{n>F(\delta)}U_n^{\rm b}(\delta)
$$
has measure
\beq \label{tail}
    \nu(U_{\omega}^{\rm b}) = o(\delta).
\eeq
In fact, the derivation of the Local Ergodic Theorem from the tail bound
does not require the Ansatz, so we will not repeat it here, see
\cite[Section~5]{KSS90}. In what follows, we prove the tail bound.

First we need a few additional constructions. Denote
\begin{align*}
   \hat{U}_n^{\rm b}(\delta) &=\big\{x\in U_0\big|\;
   \exists\text{ a divergent LOM }\Sigma,
    -T^nx\in\partial\Sigma,\; \partial\Sigma
   \cap\cS_0=\emptyset,\\ &\quad\ \ T^n \text{ and } T^{-1}
   \text{ are smooth on }\text{int}\,\Sigma,\;
   T^{-1}\Sigma\text{ is also divergent}, \\
   &\quad \ \ T^{-1}\text{ is not smooth at the endpoint }
   x'\in\partial\Sigma\text{ other than } -T^nx, \\
   &\quad \ \ \text{dist}(-T^nx,\,x')\le\kappa_{n,c_3\delta}
   (x)^{-1}c_3\delta, \text{ and }T^nx'\in U_0\big\}.
\end{align*}
For any $x\in\hat{U}_n^{\rm b}(\delta)=\hat{U}_n^{\rm b}$ and $\Sigma$ as above,
we denote
\begin{align}
&z\left(T^nx,\,\Sigma\right)=\text{dist}(-T^nx,x')\big|\;
x'\in\partial\Sigma, \; x'\ne-T^nx, \nonumber\\
& T^{-1} \text{ is not smooth at }x' \quad
\left(\le\kappa_{n,c_3\delta}(x)^{-1}c_3\delta\right). \label{2.5}
\end{align}

We note that $z_{\rm tub}(T^nx) \le z\left(T^nx,\,\Sigma\right)$.

For any point $x\in\hat{U}_n^{\rm b}$ as above, we choose a phase point
$x_{\eps_1}=-\Phi^{\eps_1}(T^nx)$ with a suitably selected $\eps_1$,
$0<\eps_1<\tau(T^nx)$.

Note. From now on we will be recycling the notation $\eps_1$ that
appeared earlier in the closed formula in Theorem 1. We think that
this action should not be the source of any confusion.

For an additional condition on how to select
$\eps_1$, see below. For any $\Sigma$ featuring the definition of
$\hat{U}_n^{\rm b}(\delta)$ and (\ref{2.5})
let $\hat{\Sigma}$ denote the flow-sync version
of $\Sigma$ containing the point $x_{\eps_1} = (q_{\eps_1}, \,
v_{\eps_1}) = -\Phi^{\eps_1} (T^nx)$, see Fig. 1.
Now $x_1 = (q_1, v_1) \in
\partial \hat{\Sigma}$ is the projection (by the flow) of the point
$x'\in\partial\Sigma$ defined above, with the property
$$
   \text{dist}(-T^nx,\,x')=z(T^n x,\,\Sigma).
$$

The other endpoint of the curve $\hat{\Sigma}$ is
$x_{\eps_1}=(q_{\eps_1},\, v_{\eps_1})$, see Fig. 1.

While selecting the time $\eps_1$ above, we try to make it sure that $x_1$ be
a post-singularity phase point, i.e. $T(-x_1)\in \cS_0$ and
$Tx_1\not\in\cS_0$.

Here we first consider the case when such a synchronization of
$\hat{\Sigma}$ is possible. After that, right before exposing (6), we
explain how to modify the following argument if the required
synchronization is not feasible.

\medskip

Consider the line segment
$$
  \cH=\left\{q\left(\Phi^t(q,\,v_{\eps_1})\right)\big|\;
(q,\,v_{\eps_1})\in \cS_0^+,\; 0<t<\tau(q,\,v_{\eps_1})\right\}
$$
in the domain $Q$, see Fig. 1. (Note that the point $q_1$ does not
belong to $\cH$, since $v_1\ne v_{\epsilon_1}$ for any strictly
convex $\hat{\Sigma}$.) The configuration
component $q_3\in\cH$ of the phase point $x_3=(q_3,\,v_{\eps_1})$ is
defined as the orthogonal projection of $q_{\eps_1}$ onto the line
$\cH$. According to this definition, the line segment
\beq \label{2.6}
   \left\{\lambda q_{\eps_1}+(1-\lambda)q_3\big|\;
   0\le\lambda\le1\right\}
\eeq
is perpendicular to $\cH$ at $q_3$.

\begin{figure}[htb]
    \centering
    \psfrag{E}{$\hat{\Sigma}$}
    \psfrag{S}{$\cS_0^+$}
    \psfrag{b}{$x_{\eps_1}=(q_{\eps_1},v_{\eps_1})$}
    \psfrag{c}{$\tilde{q}_3$}
    \psfrag{d}{$q_3$}
    \psfrag{f}{$v_{\eps_1}$}
    \psfrag{g}{$v_3$}
    \psfrag{e}{$x_1=(q_1,v_1)$}
    \includegraphics[height=2in]{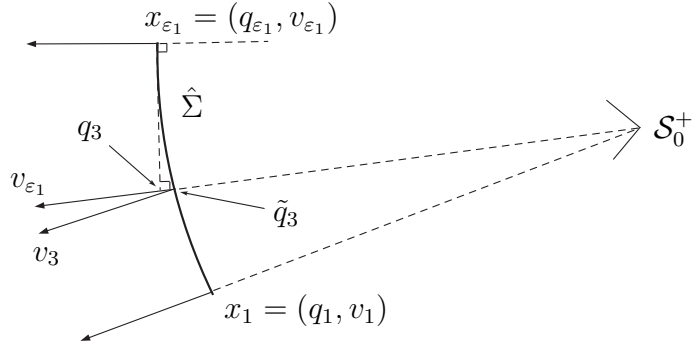}
    \caption{Illustration to Lemma~\ref{LmF}.}
    \label{Fig1}
\end{figure}

\begin{lemma} \label{LmF}
The scalar product condition of Lemma~\ref{Lm2.1} holds true for the
pair of phase points $(x_3,\,x_1)$ and $(x_3,\,x_{\eps_1})$, i. e.
$$
\aligned
&\langle q_1-q_3,\, v_1-v_{\epsilon_1}\rangle \ge 0 \\
&\langle q_{\eps_1}-q_3,\, v_{\eps_1}-v_{\epsilon_1}\rangle \ge 0,
\endaligned
$$
see also Fig.~\ref{Fig1}.
\end{lemma}

\noindent{\em Proof}. Denote the
point of intersection of the line segment $\cH$ and the carrier of
the LOM $\hat{\Sigma}$ by $\tilde{q}_3$ and the outer unit normal
vector to $\hat{\Sigma}$ at $\tilde{q}_3$ by $v_3$, see Fig. 1.
Then, by the convexity of $\hat{\Sigma}$, $q_3=\tilde
q_3+\eta v_{\epsilon_1}$ with some small scalar $\eta>0$, $\langle
q_1-\tilde q_3,\, v_3-v_{\epsilon_1}\rangle\ge 0$, and $\langle
q_1-\tilde q_3,\, v_1-v_3\rangle\ge 0$. Thus we obtain the chain of
inequalities
$$
\aligned &\langle q_1-q_3,\, v_1-v_{\epsilon_1}\rangle= \langle
q_1-\tilde q_3,\, v_1-v_{\epsilon_1}\rangle+\eta
\langle v_{\epsilon_1},\, v_{\epsilon_1}-v_1\rangle \\
&\ge \langle q_1-\tilde q_3,\, v_1-v_{\epsilon_1}\rangle = \langle
q_1-\tilde q_3,\, v_1-v_3\rangle + \langle q_1-\tilde q_3,\,
v_3-v_{\epsilon_1}\rangle \ge 0,
\endaligned
$$
finishing the proof of the lemma. \qed \medskip

Next, consider the two-dimensional ``tube'' (actually, a strip)
\beq \label{2.8}
 \cT=q\left(\left\{\Phi^tw\big|\; w\in\hat{\Sigma},\; 0\leq t
 \leq \tau_{n+1}(w)\right\}\right),
\eeq
where $q(x)$ denotes the natural projection of $\Omega$ onto $Q$, and
$\tau_{n+1}(w)$ is the time of the $(n+1)$st collision on the forward
orbit of $w$. Recall that $T^{n+1}$ is smooth on $\hat{\Sigma}$ (where
$T^{n+1}w$ is defined as $\Phi^{\tau_{n+1}(w)}(w)$), the
endpoints of $T^{n+1}(\hat{\Sigma})$ belong to the base neighborhood
$U_0$ by the construction of $\hat{\Sigma}$, and the
curves $T^{n+1}(\hat{\Sigma})$ are monotonic in the canonical $(r,\,
\phi)$ (arc-length, angle of reflection) coordinates, thus all the
``landing points'' $\Phi^{\tau_{n+1}(w)}w=T^{n+1}w$
($w\in\hat{\Sigma}$) are in the base neighborhood $U_0$, hence these
points $w$ are all sufficient.

\begin{lemma} \label{footpoint_lemma}
The footpoint $q_3=q(x_3)$ belongs to the strip $\mathcal{T}$.
\end{lemma}

\noindent{\em Proof}.
Drop a perpendicular line $l$ from $q_3$ to the supporting curve
$q(\hat{\Sigma})$ of $\hat{\Sigma}$. It follows from the convexity of
$\hat{\Sigma}$ that the intersection point $q_4$ of $l$ and
$q(\hat{\Sigma})$ lies on the arch connecting $q_{\epsilon_1}$ and
$\tilde{q}_3$, and $x_3=\Phi^\eta(q_4,v_4)$ with some small $\eta>0$
and $x_4=(q_4,v_4)\in\hat{\Sigma}$. \qed

According to Lemma~\ref{Lm2.1}, some
small time shifts of the phase points $x_3$ and
$x_1\in\partial\hat{\Sigma}$ are contained by a divergent LOM, thus
they are not mapped to the same foot point by any positive iterate
of the flow, as long as that iterate is smooth on the mentioned LOM.
(In other words, the billiard flow, being semi-dispersive, lacks
focal points.)

The same statement can be made regarding the pair of phase points
$(x_3,\,x_{\eps_1})$. Therefore, the orbit segment
$q\left(\Phi^{[0,\,\tau_{n+1}(x_3)]}x_3\right)$ of $x_3$ does not
intersect the boundary of the strip $\cT$ of (\ref{2.8}), so the
orbit segment $q\left(\Phi^{[0,\,\tau_{n+1}(x_3)]}x_3\right)$ cannot
escape from $\cT$. As a consequence, its endpoint
$\Phi^{\tau_{n+1}(x_3)} x_3=T^{n+1}x_3$ lies in $U_0$, hence the
point $x_3$ is future sufficient.

\medskip

\begin{figure}[htb]
    \centering
    \psfrag{E}{$\hat{\Sigma}$}
    \psfrag{a}{$q_3$}
    \psfrag{b}{$x_{\eps_1}=(q_{\eps_1},v_{\eps_1})$}
    \psfrag{c}{$x_3=(q_3,v_{\eps_1})$}
    \psfrag{p}{$\partial Q$}
    \psfrag{d}{$v_{\eps_1}$}
    \psfrag{g}{$q_{\ast}$}
    \psfrag{H}{$\mathcal H$}
    \psfrag{e}{$x_1=(q_1,v_1)$}
    \includegraphics[height=2in]{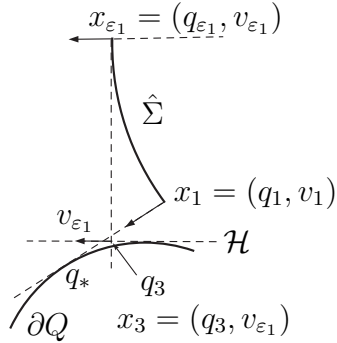}
    \caption{Illustration to the argument below.}
    \label{Fig2}
\end{figure}

As we said earlier, there is the possibility that the required
synchronization of $\hat{\Sigma}$ (with a small enough $\eps_1>0$, so
that $x_1$ becomes a post-singular phase point) is not feasible. This
phenomenon can only happen if the collision of $T^nx$ takes place very
close to a corner of the configuration space $\mathbf Q$, and very soon
after this collision the orbit segment $[T^nx,\, T^{n+1}x]$ flies near
a tangency, so that no matter how small $\eps_1>0$ one chooses, the
endpoint $x_1$ of $\hat{\Sigma}$ (other than $x_{\eps_1}$) is always a
pre-tangency phase point, see Fig. 2. In this case the proof can proceed with
such a $\hat{\Sigma}$ and $x_1$ as follows:

The forward orbit of the pre-tangency phase point $x_1=(q_1,\, v_1)$
touches $\partial\mathbf Q$ at $q_*$. Let $\mathcal H$ be the tangent
line to $\partial\mathbf Q$ near $q_*$ that is parallel to the
velocity $v_{\eps_1}$ of $x_{\eps_1}$ (analogously to the previous
case), and let again $q_3$ be the perpendicular projection of
$q_{\eps_1}$ onto $\mathcal H$, and $x_3=(q_3,\, v_{\eps_1})$, see
Fig. 2. It is clear from the picture that the forward orbit of $x_3$
enters the trip $\mathcal T$ soon after time zero. Furthermore, this
forward orbit is bound to stay in $\mathcal T$ until it reaches $U_0$
by the same reasoning as before: Both pairs $(x_3,\, x_{\eps_1})$ and
$(x_3,\, x_1)$ can be embedded in a divergent LOM by Lemma 1, thus the
footpoints $q\left(S^t x_3\right)$ of $S^tx_3$ ($t>0$, $S^tx_3$ is in
the closure of the strip $\mathcal T$) cannot be equal to any of the
footpoints $q\left(S^\tau x_{\eps_1}\right)$ or $q\left(S^\tau
x_1\right)$ ($\tau>0$), so the forward orbit of $x_3$ is unable to
escape from the strip $\mathcal T$.  Since reaching the neighborhood
$U_0$ with its forward orbit, the phase point $x_3$ proves to be
future sufficient.

\medskip

According to the previous constructions, the relevant set
for the ``upgraded'' version of Local Ergodic
Theorem is defined as  follows:

\begin{align} \label{2.9}
\tilde{U}_n^{\rm b}(\delta)
&=\Big\{x\in U_0 \big|\; \exists\eps_1,\;
0<\eps_1<\tau(T^nx),\exists
\text{ a past sufficient point} \nonumber \\
&\quad \ \ y\in\text{int}\Omega,\; Ty\in\cS_0,\;
v(y)=v(T^nx)=v\left(\Phi^{\eps_1}(T^nx)\right), \nonumber\\
&\quad\ \ \Delta q=q(y)-q\left(\Phi^{\eps_1}(T^nx)\right)\perp
v(T^nx),\; ||\Delta q||\le\kappa_{n,c_3\delta}(x)^{-1}c_3\delta\Big\}, \\
\tilde{U}_{n,m}^{\rm b}&=\tilde{U}_{n,m}^{\rm
  b}(\delta)=\left\{x\in\tilde{U}_n^{\rm b}\big|\;
\Lambda^m\le\kappa_{n,c_3\delta}(x)<\Lambda^{m+1}\right\},\nonumber
\end{align}
see the definition of the corresponding set $U_{n,m}^{\rm b}$ before Lemma
6.3 in \cite{KSS90}.

Note. For any point $x\in\hat{U}_n^{\rm b}$ the existence of a suitable
point $y$ (required by (6)) is shown by taking
$y=-x_3=(q_3,\, -v_{\eps_1})$. Hence $\hat{U}_n^{\rm
  b}\subset\tilde{U}_n^{\rm b}$.

Now, with the above constructions, we are ready to complete the
proof of the tail bound (\ref{tail}). It goes along the same lines
as in \cite[Section~6]{KSS90}, but at several points the argument
needs modifications in order to avoid using the Ansatz. We describe
these modifications in detail.

\medskip 1. First of all, in construction of local stable manifolds
\cite[Lemma~5.4]{KSS90}, we can restrict ourselves to the (limits of
the) inverse images $\Sigma_0^t(y) = \Phi^{-t} \left( \Sigma_t^t(y)
\right)$ of strictly concave, local orthogonal manifolds $\Sigma_t^t
(y)$ containing the phase point $\Phi^ty =y_t$. Indeed, if
necessary, the partially flat, concave, local orthogonal manifolds
$\Sigma_t^t(y)$ may be slightly curved to make them strictly
concave. These arbitrarily small perturbations of the manifolds
$\Sigma_t^t(y)$, obviously, produce no effect on the limiting
process and the overall proof of the Local Ergodic Theorem.

\medskip 2. The ``$n$-step bad set'' $U_n^{\rm }b=U_n^{\rm b}(\delta)$ defined above
can be replaced with our new $\hat{U}_n^{\rm b} = \hat{U}_n^{\rm b} (\delta)$, for
what it makes a phase point $y\in U_0$ ``$n$-step bad'', i.e.\
$z(T^n y)<\kappa_{n,c_3\delta}(y)^{-1}c_3\delta$, is not that this
inequality holds true, rather the fact that the construction of the
stable manifold at the phase point $y$ breaks down at the $(n+1)$-st
iteration of the billiard map, due to hitting a nearby singularity,
just as precisely described in the definition of the set
$\hat{U}_n^{\rm b}$ above. We note that, due to this change, the sets
$U^{\rm b}=\cup_{n=1}^\infty U_n^{\rm b}$ and $U^{\rm g}=U_0\setminus
U^{\rm b}$ will change, accordingly.

\medskip 3. Now, in the tail bound (\ref{tail}), the measures
of the sets $\hat{U}_n^{\rm b}$ and $\hat{U}_\omega^{\rm b} = \cup_{n>F(\delta)}
\hat{U}_n^{\rm b}$ must be estimated from above. Furthermore, these sets
may be replaced by the larger sets $\tilde{U}_n^{\rm b}$ of (\ref{2.9})
and $\tilde{U}_\omega^{\rm b} = \cup_{n>F(\delta)} \tilde{U}_n^{\rm b}$.

\medskip
4. The centerpiece estimate in the proof of the tail bound
\cite[Section~6]{KSS90} is that for any given $m\in\mathbb{N}$ the $\nu$
measure of the set
\beq \label{3.1}
  \bigcup_{n\ge N_\eta} T^n\tilde{U}^{\rm b}_{n,m}\subset \left(
  \cS_1\setminus K_\eta\right)^{[c_3\delta]}
\eeq
(featuring (6.10) in \cite{KSS90}) is bounded above by $c_\eta\delta$
if $\delta$ is small enough, where the constant $c_\eta>0$ can be made
arbitrarily small by choosing the parameter $\eta>0$ small enough.
Here the compact subset $K_\eta$ of the singularity manifold $\cS_1$
almost exhausts $\cS_1$, that is,
\beq \label{3.2}
  \bm_{\cS_1}\left(\cS_1\setminus K_\eta\right)<\eta,
\eeq
where $A^{[c_3\delta]}$ denotes the open $(c_3\delta)$-neighborhood of a
subset $A\subset\cS_1$ inside $\cM$, $\bm_{\cS_1}$ is the Lebesgue
measure on $\cS_1$, and $N_\eta \nearrow\infty$ (as $\eta\to 0$) is
some threshold function.
We note that, both in the original proof of the
Fundamental Theorem and in the current one, the set $K_\eta$ consists
of sufficient points only, and -- in the original proof -- exhausting
$\cS_1$ in the sense of (\ref{3.2}) was made possible by the Ansatz!
On the other hand, in the current scenario, lacking the Ansatz, we can
only say that $K_\eta \subset \text{Suff} \left(\cS_1\right)$, where
$\text{Suff}\left(\cS_1\right)$ denotes the set of all sufficient
points of $\cS_1$, and
$$
\bm_{\cS_1}\left(\text{Suff}\left(\cS_1\right)\setminus
K_\eta\right) <\eta.
$$
We denote the set $\text{Suff}\left(\cS_1\right)\setminus K_\eta$ by
$B_\eta$.

\medskip
5. Keeping in mind definition (\ref{2.9}) of the sets
$\tilde{U}_n^{\rm b}$, some open, tubular neighborhood $V_0$ of $\cS_1$ in
$\cM$ has a uniquely defined foliation
$$
  \Psi\colon \cS_1\times(-\eps_0,\,\eps_0)
  \buildrel\cong\over\longrightarrow V_0
$$
with the properties that the section $\Psi \left( \cS_1 \times \{0\}
\right)$ is equal to $\cS_1$ (more precisely, $\Psi(y,0)=y$ for all
$y\in\cS_1$), and on the one-dimensional foliae (curves) $\Psi
\left( \{y_0\} \times (-\eps_0,\, \eps_0)\right)$ (here $y_0 \in
\cS_1$) the phase points have a constant velocity vector
\beq \label{3.5}
  v\left(\Psi(y_0,s)\right)=
  v\left(\Psi(y_0,0)\right)=
  v(y_0)
\eeq
for $|s|<\eps_0$, and the curve $\Psi(y_0,s)$ has a time-sync version
$$
\gamma_0(s)=\Phi^{\tau(s)}\left(\Psi(y_0,s)\right)
$$
(where $0< \tau (s) < \tau \left( \Psi (y_0,s) \right)$ and
$|s|<\eps_0$) such that the carrier $q(\gamma_0(s)) \subset Q$ is
linear and orthogonal to the flow-invariant hull
$$
 \left\{q\left(\Phi^\tau(y_0)\right)\big|\;
 0<\tau<\tau(y_0)\right\},
$$
and $\left\Vert\dfrac{d}{ds}q\left(\gamma_0(s)\right)\right\Vert=1$
for all $y_0\in\mathcal{S}_1$ and $|s|<\epsilon_0$. (So $s$ is an
arc-length parameter.)  Furthermore, according to (\ref{2.9}), the
crucial set
$$
   \cup_{N\ge N_\eta}T^n\tilde{U}^{\rm b}_{n,m}
$$
in (\ref{3.1}) is a subset of
$$
  \Psi\left (B_\eta\times[-c_3\delta,\,c_3\delta]\right)
$$
with $\bm_{\cS_1}(B_\eta)<\eta$, so we have that
$$
  \nu\left(\cup_{N\ge N_\eta}T^n\tilde{U}^{\rm b}_{n,m}\right)\leq
  c_2c_3\eta\delta
$$
with an absolute constant $c_2>0$. The existence of such a constant
$c_2$ is attributed to the facts that

\medskip

(a) the flow invariant measure $\mu$ is uniform on $\Omega$, i.e.\
$d\mu=\text{const}\cdot dq \, d\theta$, where $\theta$ is the
angular coordinate of the unit velocity vector;

\medskip

(b) the $T$-invariant measure $\nu$ on $\cM$ is the projection of
$\mu$ onto $\cM$ via the flow;

\medskip

(c) the above curves $\gamma_0(t)$ are uniformly transversal to the
forward invariant hull $\cup_{t>0} \Phi^t (\cS_0)$ of the
singularity manifold $\cS_0$.

\medskip

As a matter of fact, this constant $c_2$ is exactly the same as the
constant appearing in \cite[Lemma~2]{SC87}, and the proof of its
existence goes along the same lines, see also
\cite[Lemma~4.10]{KSS90}.

Since the multiplier $c_2c_3\eta$ of $\delta$ can be made
arbitrarily small by choosing the number $\eta>0$ small enough, we
finish the proof of the tail bound without having used the Ansatz.

\end{document}